\documentclass[12pt]{amsart}

\usepackage[leqno]{amsmath}
\usepackage[all]{xy}
\usepackage[mathscr]{euscript}
\usepackage{amssymb}
\usepackage{amscd}
\usepackage{amsthm}
\usepackage{epsfig}
\usepackage{a4wide}

\usepackage{graphicx}

\usepackage[colorlinks=true,citecolor=black,linkcolor=black,urlcolor=blue]{hyperref}

\theoremstyle{plain}
\newtheorem{theorem}{Theorem}
\newtheorem{lemma}[theorem]{Lemma}
\newtheorem{corollary}[theorem]{Corollary}
\newtheorem{proposition}[theorem]{Proposition}

\theoremstyle{definition}

\theoremstyle{remark}

\newcommand{\pf}{\noindent{\em Proof: }}
\newcommand{\epf}{\hfill\hbox{\rule{3pt}{6pt}}\\}

\newcommand{\R}{{\mathbb R}}
\newcommand{\G}{{\mathcal G}}
\newcommand{\N}{{\mathbb N}}
\newcommand{\s}{\mathcal S}
\newcommand{\I}{\mathcal I}

\newcommand{\B}{\mathcal B}

\title[The polytopal structure of the tight-span]{\bf The polytopal structure of the tight-span of a totally split 
decomposable metric}

\author{K.~T.~Huber}
\address[K.~T.~Huber]{School of Computing Sciences, University of East Anglia,
United Kingdom.}
\email{k.huber@uea.ac.uk}

\author{J.~Koolen}
\address[J.~Koolen]{Wen-Tsun Wu Key Laboratory of CAS,
  School of Mathematical Sciences, University of
  Science and Technology of China, Hefei, Anhui, 230026, P.R. China.}
\email{koolen@ustc.edu.cn}

\author{V.~Moulton}
\address[V.~Moulton -- Corresponding author]{School of Computing Sciences, University of East Anglia,
United Kingdom}
\email{v.moulton@uea.ac.uk}

\date{\today}

\begin{document}

\maketitle


\begin{abstract}
  The tight-span of a finite metric space is a polytopal complex that 
has appeared in several areas of mathematics. In this paper we
determine the polytopal structure of the tight-span
of a totally split-decomposable (finite) metric. Totally split-decomposable
metrics  are a generalization of tree-metrics and have 
importance within phylogenetics.
In previous work, we showed that the cells of the tight-span of 
such a metric are zonotopes that are polytope isomorphic to 
either hypercubes or rhombic dodecahedra.
Here, we extend these results and show that the tight-span
of a totally split-decomposable metric can be
broken up into a canonical 
collection of polytopal complexes whose polytopal structures can 
be directly determined from the metric. 
This allows us to also completely determine the
polytopal structure of the tight-span of a totally 
split-decomposable metric in a very direct way.
We anticipate that  our improved understanding of this structure 
may ultimately lead to improved techniques for phylogenetic inference.

\end{abstract}

\section{Introduction}

In this paper, $X$ will denote a finite set with $|X|\ge2$. Given 
a metric $d$ on $X$, the {\em tight-span} of $d$ is the polytopal complex
$T(d)$ which consists of the bounded faces of the polyhedron
$$
\{ f \in \R^X \,: \, f(x) + f(y) \ge d(x,y), \mbox{ for all } x,y \in X \}.
$$ 
The tight-span of an arbitrary metric was first introduced by 
Isbell \cite{I64} (where it was called the {\em injective hull}), and was subsequently 
redisovered in \cite{CL94,D84}. It has appeared in 
various areas of mathematics 
including group theory \cite{D89,L13}, phylogenetics \cite{BD92,DHM02b}, network 
flow theory \cite{H09}, tropical geometry \cite{sturmfels} and the theory of
low distortion embeddings \cite{BC96}. 

\begin{figure}[h]
	\begin{center}
		\includegraphics[width=10cm]{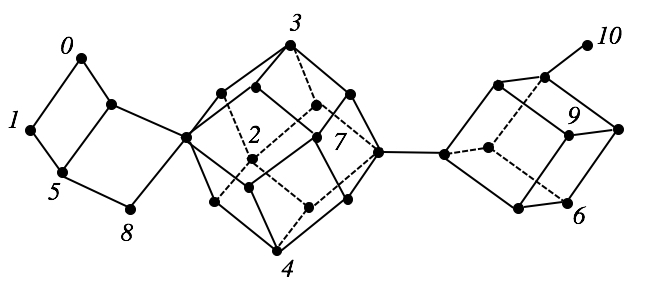}
	\end{center}
	\caption{An example of the 1-skeleton of the tight-span of a totally 
split-decomposable metric defined on the set $\{0,1,2,\dots,9,10\}$ (the
metric is defined in the next section). It is a 
polytopal complex
		consisting of 5 blocks: two 1-cubes, a 3-cube, a rhombic 
		dodecahedron and a pair of 2-cubes with an edge in common. }\label{tspan}
\end{figure}

In this paper we are interested in determining the 
polytopal structure of the tight-span of a special type
of metric on a finite set called a  {\em totally split-decomposable} metric \cite{BD92}.
We picture the 1-skeleton for an  example of the
tight-span of such a metric in Figure~\ref{tspan}.
This type of metric originated from the study of tree-metrics arising 
in evolutionary biology, and is now commonly used in the phylogenetic 
analysis of molecular DNA sequence data \cite{huson}. Recently, infinite 
versions of these metrics have also appeared in the study 
of certain finitely generated groups \cite{P15}.

The tight-span of a totally split-decomposable metric has 
an interesting polytopal 
structure.  Some of the first results concerning this structure appeared
in \cite{DHM01b,DHM02b} where, amongst other things,  
it was shown that the tight-span of a certain subclass
of the totally split-decomposable metrics  can be considered as a special 
type of median complex called the {\em Buneman complex}
(cf. \cite{C00} for more concerning median complexes). 
Subsequently, in \cite{HKM06} it was shown that the cells in 
the tight-span (i.e. the polytopes from which the tight-span is comprised) 
of a totally split-decomposable metric are
zonotopes that are polytope isomorphic to either hypercubes or 
rhombic dodecahedra. In this paper, we shall extend
this result by completely determining the polytopal structure
of the tight-span of a totally split-decomposable metric.

We now summarize the rest of the paper. In 
Section~\ref{sec:incomp} we begin by presenting some
terminology and results concerning split systems, 
structures which form the basis for defining totally split-decomposable metrics.
We then present a structural result (cf. Corollary~\ref{structure}) concerning 
a certain graph that can be associated to a
split system 
(called the {\em incompatibility graph}). This 
result provides a key to breaking up the tight-span 
of a totally split-decomposable metric into easier to
understand pieces.

In Section~\ref{sec:buncomp},
we study properties of the Buneman complex, a polytopal 
complex that can be associated to a weighted split system $(\s,\alpha)$.
As we shall explain, we can break this complex down into pieces or blocks
(maximal connected subcollections of cells, which cannot
be disconnected by removing a vertex or 0-cell)
which are in bijective correspondence with the 
connected components of the incompatibility graph of $\s$.
In Section~\ref{sec:key}, we then prove a result concerning 
the blocks of the Buneman complex (Theorem~\ref{pairs}), which 
ultimately allows us in Section~\ref{sec:blocks} 
to reduce the problem of understanding the polytopal 
structure of the tight-span of a totally split-decomposable 
metric to that of understanding its blocks.

In Section~\ref{sec:kappa} we begin to relate properties of the Buneman complex 
to those of the tight-span. We do
this by considering a map $\kappa$ that maps the Buneman complex onto
the tight-span, which was defined in \cite{DHM98} and 
further studied in  \cite{HKM06}. In particular, we show that 
the map $\kappa$ has certain properties relative to 
the blocks of the Buneman complex, which allows
us to prove in Section~\ref{sec:relate} that the map $\kappa$
induces a bijection between the blocks of the Buneman complex
and the blocks of the tight-span (Theorem~\ref{Kmap}). 

In Section~\ref{sec:blocks}
we conclude by showing that, for most of the 
blocks in the Buneman complex, $\kappa$ induces a polytopal complex isomorphism
between the blocks in the Buneman complex
and their corresponding blocks in the tight-span 
(Theorem~\ref{theo:structure}).
Moreover, we see that the remaining blocks in the tight-span have 
a very simple structure:
they are all rhombic dodecahedra (Theorem~\ref{Kmap}). As we shall also explain
in Section~\ref{sec:blocks}, this allows us to completely
determine the polytopal structure of the tight-span of a totally
split-decomposable metric.

In future work it could be of interest to understand how our results
may extend to the case of infinite, totally split-decomposable
metrics defined in \cite{P15}. In addition, it could be of interest 
to understand how tight-spans of totally split-decomposable metrics
fit the theory of CAT(0) complexes \cite{C00}. 
Ultimately, we anticipate that better understanding the structure
of the tight-span of a totally split-decomposable metric might be 
useful within phylogenetics. For example, this structure 
is closely-related to the block-realization of a metric \cite{DHKM08},
which may give insights on how to decompose metrics
that are of importance in phylogenetics.

Throughout the paper we will follow and extend the 
notation and definitions presented in \cite{HKM06}. For 
the reader's convenience, we shall briefly recall relevant notation from \cite{HKM06},
but we refer the reader to that paper for more detail.

\section{The incompatibility graph of weakly compatible split system}\label{sec:incomp}

\begin{figure}[h]
	\begin{center}
		\includegraphics[width=10cm]{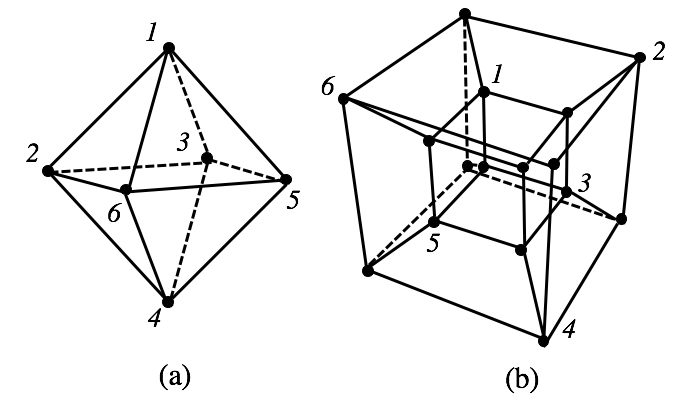}
	\end{center}
	\caption{(a) An octahedral split system on the set $\{1,2,3,4,5,6\}$ 
		consisting of the 
		splits 
		$\{1,2,3\} | \{4,5,6\}$,
	    $\{2,3,4\}| \{5,6,1\}$, 
	    $\{3,4,5\}|\{6,1,2\}$, and 
	    $\{1,3,5\}|\{2,4,6\}$, 
        in which each split is obtained by taking the labels of some
face and its opposite face in the pictured octahedron. (b) 
The Buneman complex associated to the split system in (a).}\label{os}
\end{figure}

We begin by recalling some terminology and results concerning split systems,
structures which form the basis for defining and understanding totally split-decomposable
metrics. A {\em split} of $X$ is a  bipartition of $X$, and a set $\s$ of splits of 
$X$ is a {\em split system (on $X$)}. We denote a split $\{A,B\}$
of $X$ with $\emptyset\not=A,B\subset X$ by $A|B$ ($=B|A$).
For a split $S$ of $X$ and some elements
$x\in X$ we denote by $S(x)$ the element of $S$ that contains $x$
and by 	$\overline{S}(x)$ the complement of $S(x)$ in $X$.
Two distinct splits $S,S' \in  \s$ are {\em compatible} if there
exists $A \in S$ and $A' \in S'$ such that $A \cup A' = X$, otherwise
$S$ and $S'$ are {\em incompatible}. We call a split system $\s$ {\em incompatible}
if every pair of distinct splits in $\s$ is incompatible. We also 
define a split system with 1 element to be incompatible.
Following \cite{BD92}, we call  a split system {\em weakly compatible} if there exist no three splits 
$S_1,S_2,S_3 \in \s$ and four elements $x_0,x_1,x_2,x_3 \in X$ 
such that 
\begin{eqnarray}\label{wc-condition}
S_j(x_i) = S_j(x_0) \mbox{ if and only if } i=j.
\end{eqnarray}
Note that in \cite{DHM00} weakly compatible split systems
where characterized as those split system $\s$ for which
for any three splits $S_1,S_2,S_3\in\s$ and all $x\in X$, we have
$\overline{S_1}(x)\cap \overline{S_2}(x) \cap \overline{S_3}(x)
\in \{ \overline{S_1}(x)\cap \overline{S_2}(x),
	\overline{S_2}(x)\cap \overline{S_3}(x),
	\overline{S_1}(x)\cap \overline{S_3}(x)  \}.$
Also note that			
in \cite[Theorem 4.1]{DHM00},  the following result 
is proven. Suppose that $\s=\{S_1,\ldots, S_k\}$ 
is a weakly compatible yet incompatible split system
of size $1\leq k$. Then there exists a partition 
$X=X_1\stackrel{\cdot}{\cup}\ldots \stackrel{\cdot}{\cup} X_{2k}$ 
of $X$ into $2k$ non-empty pairwise disjoint subsets $X_i$
such that either $\s$ is {\em strictly circular}, that is
$S_i=X_i
\stackrel{\cdot}{\cup}\ldots \stackrel{\cdot}{\cup} X_{i+k}|
X_{i+k+1}\stackrel{\cdot}{\cup}\ldots \stackrel{\cdot}{\cup} X_{i-1}
$
holds for $1\leq i\leq k$ or
$\s$ is {\em octahedral}, that is, $k=4$ and we can relabel the
elements in $\s$ such that
$S_i=X_i\stackrel{\cdot}{\cup} X_{i+1}\stackrel{\cdot}{\cup} X_{i+2}|
X_{i+3}\stackrel{\cdot}{\cup} X_{i+4}\stackrel{\cdot}{\cup} X_{i+5}$
for $1\leq i\leq 3$  and 
$S_4=X_1\stackrel{\cdot}{\cup} X_3\stackrel{\cdot}{\cup} X_5|
X_2\stackrel{\cdot}{\cup} X_4\stackrel{\cdot}{\cup} X_6$ (where we take indices modulo $2k$).
See Figure~\ref{os} for a diagrammatic representation of 
an  octahedral split system. For $\s$ a split system on $X$, let 
$Oct(\s) = \{ \s' \subseteq \s \,:\, \s' \mbox{ is octahedral}\}$. 
We call $\s$ {\em consistent} if $\s$ is weakly compatible
and does not contain an octahedral subsystem, that is, if $Oct(\s)$ is empty.

Using this last result we now present a key property of 
weakly compatible split systems.

\begin{theorem}\label{intersection-result}
Suppose that $\s\subseteq \s(X)$ is a weakly compatible 
split system and that $\s_1\subseteq \s$
is an octahedral split system. Then every split $S$ in $\s-\s_1$ is
compatible with every element in $\s_1$. 
\end{theorem}
\pf Suppose $\s_1=\{S_1,S_2,S_3,S_4\}$ 
and assume for contradiction that there exists
some split $S_5\in \s-\s_1$ such that
$S_5$ is incompatible with some split in $\s_1$. 
Then 
$$
\s'=\{S\in \s_1 : S_5 \mbox{ and } S \mbox{ are incompatible }\}\not=\emptyset.
$$
Assume without loss of generality that 
$S_i=X_i\dot{\cup}X_{i+1}\dot{\cup}X_{i+2}| 
X_{i+3}\dot{\cup}X_{i+4}\dot{\cup}X_{i+5}$, $i=1,2,3$
and $S_4=X_1\dot{\cup}X_3\dot{\cup}X_5|X_2\dot{\cup}X_4\dot{\cup}X_6$
for  $X=X_1\dot{\cup}\ldots\dot{\cup} X_6$.

If $|\s'|=1$ then we may assume without loss of generality that
$\s_1=\{S_1\}$. But then $S_5$ is compatible with all $S_i$, $i=2,3,4$.
Since $\{S_2,S_3,S_4\}$ is strictly circular, there exists some
$x\in X$ such that
$S_5(x)\subsetneq S_i(x)$, $i=2,3,4$. Hence,
$S_5(x)\subseteq \bigcap_{i=2,3,4}S_i(x)=X_j$, for some $j\in \{1,\ldots,6\}$. 
Thus, $S_5$ and $S_1$ cannot be incompatible since either
$X_j\subseteq S_1(x)$ or $X_j\subseteq \overline{S_1}(x)$ 
holds; a contradiction.\\

Now suppose  $|\s'|\geq 2$ and let $x\in X$. Without loss of
generality assume that $x\in X_1$.
If $|\s'|=2$ then we may assume without loss of generality that
$\s_1=\{S_1,S_2\}$. Hence,
\begin{eqnarray}\label{intersection}
\overline{S_1}(x)\cap \overline{S_2}(x)=
\overline{S_3}(x)\cap \overline{S_4}(x).
\end{eqnarray}
Since $\{S_1,S_2, S_5\}$ is weakly compatible and
incompatible, there exist $i,j\in\{1,2,5\}$ such that
\begin{eqnarray}\label{wc-and-empty}
\,\,\,\,\,\,\,\,\,\,\emptyset\not=\overline{S_i}(x)\cap \overline{S_j}(x)
=
\overline{S_1}(x)\cap \overline{S_2}(x)\cap \overline{S_5}(x)=
\overline{S_3}(x)\cap \overline{S_4}(x)\cap \overline{S_5}(x).
\end{eqnarray}
and so $\overline{S_3}(x)\cap \overline{S_5}(x)\not=\emptyset$
 and $\overline{S_4}(x)\cap \overline{S_5}(x)\not=\emptyset$.
Since $S_5$ and $S_3$ are compatible it follows that 
either $\overline{S_5}(x)\cap S_3(x)=\emptyset$ or 
$\overline{S_3}(x)\cap S_5(x)=\emptyset$ must hold.  If 
$\overline{S_5}(x)\cap S_3(x)=\emptyset$ then
$\overline{S_5}(x)\cap S_4(x)\not=\emptyset$ as otherwise
(\ref{intersection}) would imply 
$\overline{S_5}(x)\subseteq \overline{S_1}(x)$ which is impossible.
Since $S_5$ and $S_4$ are compatible $\overline{S_4}(x)\cap S_5(x)=\emptyset$
follows. Hence, $\overline{S_4}(x)\subseteq  \overline{S_5}(x)$.
Combined with $\overline{S_5}(x)\cap S_3(x)=\emptyset$, it follows that
$\overline{S_4}(x)\subseteq \overline{S_3}(x) $ which is impossible.
Hence, $\overline{S_3}(x)\cap S_5(x)=\emptyset$. 

Since $S_5$ and $S_4$ 
are compatible either $\overline{S_5}(x)\cap S_4(x)=\emptyset$ or
$\overline{S_4}(x)\cap S_5(x)=\emptyset$  
must hold. If $\overline{S_5}(x)\cap S_4(x)=\emptyset$ held, then
$\overline{S_3}(x)\subsetneq \overline{S_5}(x)\subsetneq \overline{S_4}(x)$
which is impossible. Thus $\overline{S_4}(x)\cap S_5(x)=\emptyset$. 
But then $S_5(x)\subseteq S_4(x)\cap S_3(x)$.
Since $\s_1$ is weakly compatible yet incompatible (\ref{intersection})
implies 
$$
\bigcap_{i=1,\ldots,4}\overline{S_i}(x)
=\overline{S_1}(x)\cap \overline{S_2}(x) 
=\overline{S_3}(x)\cap \overline{S_4}(x). 
$$ 
Hence, by
\cite[Lemma 2.1]{DHM00} 
$$
\bigcap_{i=1,\ldots,4}S_i(x)
=S_1(x)\cap S_2(x) 
=S_3(x)\cap S_4(x). 
$$ 
Consequently, $S_5(x)\subseteq S_1(x)\cap S_2(x) $ and so $S_5(x)\subseteq S_1(x)$
which is impossible.\\

If $|\s'|=3$, then we may assume without loss of generality that
$\s'=\{S_1,S_2,S_3\}$. But then  we have again
\begin{eqnarray*}
\overline{S_1}(x)\cap \overline{S_2}(x)=
\overline{S_3}(x)\cap \overline{S_4}(x)
\end{eqnarray*}
and similar arguments as in the previous case imply that 
$\overline{S_4}(x)\cap \overline{S_5}(x)\not=\emptyset$.
Since $S_4$ and $S_5$ are compatible either 
$\overline{S_4}(x)\cap S_5(x)=\emptyset$ or 
$\overline{S_5}(x)\cap S_4(x)=\emptyset$  must hold. 
If $\overline{S_4}(x)\cap S_5(x)=\emptyset$ then
 we have 
 $X_2\dot{\cup}X_4\dot{\cup}X_6\subsetneq \overline{S_5}(x) $.
Hence, there exists some 
$y\in \overline{S_5}(x)-(X_2\dot{\cup}X_4\dot{\cup}X_6)$.
Since $X=X_1\dot{\cup}\ldots \dot{\cup}X_6$ it follows that there must exist
some $i\in\{1,3,5\}$ with $y\in X_i$. Without loss of generality we may assume
that $i=1$. Then the four elements
$y$, $x_3\in X_3$, $x_4\in X_4$, and $x_5\in X_5$ together with the three 
splits $S_1,S_2,S_5\in \s$ violate Property~(\ref{wc-condition}); 
a contradiction. Hence,   
$\overline{S_5}(x)\cap S_4(x)=\emptyset$ must hold. Thus
$S_4(x)\subsetneq S_5(x)$ and so 
$X_1\dot{\cup}X_3\dot{\cup}X_5\subsetneq S_5(x) $.
Replacing $\overline{S_5}(x)$ by $S_5(x)$
and $X_2\dot{\cup}X_4\dot{\cup}X_6$ by $X_1\dot{\cup}X_3\dot{\cup}X_5$
in the previous argument shows that we can again find 4 
elements in $X$ and three splits in $\s$ such that Property~(\ref{wc-condition}) is violated.\\

If $|\s'|=4$, then $\s'=\s_1$. Hence, $\s'\cup\{S_5\}$ is a weakly compatible
yet incompatible split system that contains an octahedral split system.
By \cite[Theorem 4.1]{DHM00}, it follows that $\s'\cup\{S_5\}$
is octahedral. But then $|\s'\cup\{S_5\}|=4$ which is impossible.
\epf

\begin{figure}[h]
	\begin{center}
		\includegraphics[width=10cm]{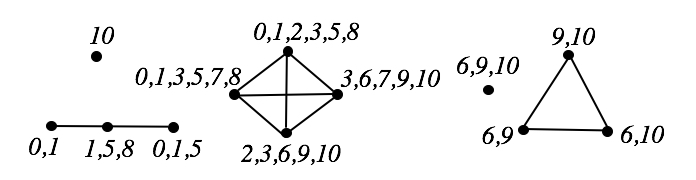}
	\end{center}
	\caption{The incompatibility graph of a split system 
on the set $\{0,1,\dots,9,10\}$ underpinning the tight-span 
pictured in Figure~\ref{tspan}. Each vertex corresponds to the split 
$A| \overline A$, with $A$ being the set of numbers labelling the vertex 
(so for, example, the vertex labelled $6,9,10$ corresponds to the 
split given by taking $\{6,9,10\}$ and its complement).}\label{ic}
\end{figure}

Given a split system $\s$ we define the {\em incompatibility graph} $
{\mathcal I}(\s)$ associated to $\s$
to be the graph with vertex set $\s$ and edge set consisting of those
pairs $\{S,S'\}$ of distinct splits $S,S' \in \s$ which are incompatible 
(cf. e.g. \cite{DHKM11}).
We also let $C(\mathcal{I}(\s))$ denote the set of connected components 
of $\mathcal{I}(\s)$. To illustrate these definitions, let  $\s$ denote
the split system on $X=\{0,1,\ldots, 9,10\}$
underpinning the totally split-decomposable metric whose
tight-span we picture in Figure~\ref{tspan}. Then $\mathcal I(\s)$ is
the graph depicted in Figure~\ref{ic}.

\begin{corollary}\label{structure}
Suppose that $\s$ is a weakly compatible split system, $\s'\subseteq \s$ is an 
octahedral split system, and ${\mathcal I}(\s)$ is connected. Then $\s=\s'$. 
In particular, if $\s' \in C({\mathcal I}(\s))$, then 
$\s'$ is octahedral or consistent.
\end{corollary}
\pf Suppose to the contrary that $\s\not=\s'$. Let $S$ denote a split in 
$\s-\s'$ and let $S'$ denote a split in $\s'$. Since, by assumption, 
$\I(\s)$ is connected, there exists some path $P$ in ${\mathcal I}(\s)$ from $S$ to $S'$.
But then there must exist some edge $\{S_1,S_2 \}$ on $P$ such that 
$S_1\in \s-\s'$ and $S_2\in \s'$. But this is impossible
since then $S_1$ and $S_2$ are incompatible in contradiction to
Theorem \ref{intersection-result}.
\epf

\section{The Buneman complex of a weighted split system}\label{sec:buncomp}

In this section we gather together some useful results concerning the 
Buneman complex. We begin by presenting some general definitions 
concerning polyhedral complexes (cf. e.g. \cite{HB}).

\subsection{Polytopal complexes}

A {\em polyhedron} in $\R^n$, $n \in \N$, is the intersection of a 
finite collection of halfspaces in $\R^n$ and a  {\em polytope} is
a bounded polyhedron. A {\em face} of a polyhedron $P$ is
the empty-set, $P$ itself, or the intersection of $P$ with a
supporting hyperplane. We denote the fact that $P$ is 
$d$-dimensional by putting $dim(P)=d$. The 0-dimensional 
faces of $P$ are also called its {\em vertices} and the 1-dimensional 
faces its {\em edges}.
A {\em polyhedral complex}
$\mathcal C$ is a finite collection of polyhedra (which we call {\em cells}) 
such that each face of a member of $\mathcal{C}$  is itself a member
of $\mathcal{C}$, and the intersection of two members of $\mathcal{C}$ 
is a face of each. If all members of $\mathcal{C}$  are polytopes, 
we call $\mathcal{C}$ a {\em polytopal complex}.
The  {\em 1-skeleton }  of $\mathcal{C}$ is the union
of its 0- and 1-dimensional cells (which we will also 
consider sometimes as being a graph, whose vertices
and edges correspond to vertices and edges in $\mathcal C$). 
Note that we will not usually distinguish between $\mathcal{C}$ 
and its underlying set $\bigcup_{C \in \mathcal{C}} C$.
For any $c$ in the underlying set of $\mathcal{C}$, we let $[c]$ denote the 
minimal cell $C$ in $\mathcal{C}$ (under cell inclusion) that contains $c$.
In this case we also call $c$ a {\em generator} of $C$.

Suppose that $\mathcal{C}$ is a connected polytopal complex.
A vertex in $\mathcal{C}$ is a {\em cut-vertex} if $\mathcal{C} - \{v\}$
is disconnected. Note that if $v$ is a cut-vertex of $\mathcal{C}$ and $C$ is some connected
component of $\mathcal{C} - \{v\}$, then $C^{+v} = C \cup \{v\}$  
can be regarded as a connected polytopal complex in the obvious way.
A maximal collection of cells in $\mathcal{C}$ that is connected and does not
contain a cut-vertex is called a {\em block} of $\mathcal{C}$. We
denote the set of blocks of $\mathcal C$ by $\B(\mathcal C)$.
Note that this is also a polytopal complex.

\subsection{The Buneman complex}

A {\em weighted split system $(\s,\alpha)$ (on $X$)}, is a split system $\s$ on $X$
together with a map $\alpha:\s \to \R^{>0}$.  
We now define the 
Buneman complex of such a split system (cf. \cite[Section 2.3]{HKM06}). First, we put
$$
\mathcal U(\s) = \{ A \subseteq X \,: \mbox{ there exists $S \in \s$ with $A \in S$ } \}.
$$
Also, given any map $\phi: \mathcal U(\s) \to \R$, we define
$$
supp(\phi) = \{A \in \mathcal U(\s) \,:\, \phi(A) \not= 0\}
$$
and put 
$$
\s(\phi) = \{ S \in \s \,:\, S \subseteq supp(\phi) \}.
$$
Now, we let 
$$
H(\s,\alpha) = \{\phi \in \R^{\mathcal U(\s)} \,:\, \phi(A) \ge 0 
\mbox{ and }
\phi(A) + \phi(\overline{A}) = \frac{\alpha(A|\overline{A})}{2}
 \mbox{ for all } A \in  \mathcal U(\s)\}
$$
(a polytope in $\R^{\mathcal U(\s)}$ which is polytope isomorphic to an $|\s|$-dimensional hypercube)
and define the {\em Buneman complex  associated to $(\s,\alpha)$ } to be the polytopal complex
$$
B(\s,\alpha) = \{ \phi \in H(\s,\alpha) \,:\, A_1, A_2 \in supp(\phi) \mbox{ and } A_1 \cup A_2 = X
\Rightarrow A_1 \cap A_2 = \emptyset \}.
$$
Note that $X$ can be considered as a subset of the set of 
vertices of $B(\s,\alpha)$ via the mapping which takes each $x \in X$ to the map
$$
\phi_x: \mathcal U(\s)\to \mathbb R^{\geq 0}
$$ 
given by $\phi_x(A)=\alpha(A|\overline A)/2$ if $x\not\in A$
and $0$ otherwise, for all $A\in\mathcal U(\s)$.
For example, for the split system in Figure~\ref{os}(a) on $X=\{1,2,\dots,6\}$ in which all 
splits are given weight 1, the Buneman complex is the 4-cube 
in Figure~\ref{os}(b),
where the labelled vertices correspond to the elements in $X$.

We now present two simple 
observations concerning the Buneman complex that will be useful later on.

\begin{lemma}\label{restrict}
	Suppose that $(\s,\alpha)$ is a weighted split system on $X$, 
	$\phi \in B(\s,\alpha)$ and $dim([\phi]) > 0$. Then for all 
	$S \not\in \s(\phi)$ and all vertices $\tilde{\phi}$ of $[\phi]$, 
we have $\tilde{\phi}|_S = \phi|_S$.
\end{lemma}
\pf
This follows from the definition of $B(\s,\alpha)$ (see also  
\cite[Lemma 3.1(i)]{HKM06}).
\epf

\begin{lemma} \label{incomp}
	Suppose that $(\s,\alpha)$ is a weighted split system on $X$, and
	$\phi \in B(\s,\alpha)$. Then
	any pair of distinct splits in $\s(\phi)$ is incompatible.
\end{lemma}
\pf
If not, then there exist $S_1 =A_1|B_1 \neq S_2=A_2|B_2$ 
in $\s(\phi)$ with $A_1 \cap A_2 = \emptyset$, say. So
$\overline{A_1} \cup \overline{A_2} = X$. Since  
$\overline{A_1}, \overline{A_2} \in supp(\phi)$, we obtain 
$\overline{A_1} \cap \overline{A_2} = \emptyset$ as $\phi\in B(\s,\alpha)$. 
Hence, 
$A_2 =\overline{A_1}$, which contradicts $S_1 \neq S_2$.
\epf

\subsection{The blocks of the Buneman complex}

We now study the set $\B(\s,\alpha)$ of blocks
of the Buneman complex. 
In what follows, we consider the 1-skeleton $G(\s,\alpha)$ 
of $B(\s,\alpha)$ also as being 
a weighted graph (where, in case an edge corresponds to 
a split $S \in \s$ in the natural way, it is weighted by $\alpha(S)$).
This graph is also known as the {\em (weighted) Buneman graph}  
of $(\s,\alpha)$ (cf. e.g. \cite{DHM97}), and we shall exploit 
some of its well-known properties.

In particular, note that each block of the 
Buneman complex $B(\s,\alpha)$ corresponds to the 
union of all of those cells in $B(\s,\alpha)$ 
whose 1-skeleta are contained in some
block (i.e. maximal 2-connected component) 
of $G(\s,\alpha)$. Hence
by \cite[Theorem~5.1]{DHKM11} (which gives a 1-1 correspondence 
between the blocks of the Buneman graph of a split system $\s$
and the set $C(\mathcal{I}(\s))$ of connected components 
of $\mathcal{I}(\s))$, each block in
$\B(\s,\alpha)$  corresponds to
precisely one element  $\s' \in C(\mathcal{I}(\s))$ and, in this 
case, is isomorphic as a polytopal 
complex to $B(\s',\alpha|_{\s'})$. 
We shall denote the block of $B(\s,\alpha)$ corresponding to 
$\s' \in C(\mathcal{I}(\s))$
(considered as a subpolytopal complex of $B(\s,\alpha)$) by
$B_{\s'}(\s,\alpha)$. In particular,
it follows that 
$$
\B(\s,\alpha) = 
\{B_{\s'}(\s,\alpha) \,:\, \s' \in  C({\mathcal I}(\s))\}.
$$
Associating to any cell $C$ of $B(\s,\alpha)$ the split
system $\s(C)$ induced by $C$ by deleting parallel edges of $C$ and
to any 
collection $\mathcal C$ of cells in $B(\s,\alpha)$ the split system $\s(\mathcal C) = \bigcup_{C \in \mathcal C} \s(C)$, we obtain $\s(B_{\s'}(\s,\alpha)) = \s'$.

We close this section by presenting a connection 
between the blocks of the Buneman complex of a weighted
split system $(\s,\alpha)$ 
and the incompatibility graph of $\s$.

\begin{lemma} \label{zero-one}
	Suppose that $(\s,\alpha)$ is a weighted split system on $X$
	and $\s' \in C(\mathcal I (\s))$.\\
	(i) If $\phi \in  B_{\s'}(\s,\alpha)$
	and $S \in \s - \s'$, then for all $A \in S$, 
$\phi(A) \in \{0,\alpha(S)/2\}$.\\
	(ii) If $\phi_1 \in B_{\s'}(\s,\alpha)$
	and $\phi_2 \in B(\s,\alpha) - B_{\s'}(\s,\alpha)$, then 
	there exists a split $S \in \s - \s'$ and some $A\in S$ such that 
	$\phi_1(A)=0$ and $\phi_2(A) \neq 0$.
\end{lemma}
\pf
\noindent (i): Since $\phi$ must be contained in a cell of $B_{\s'}(\s,\alpha)$
this is an immediate consequence of \cite[Lemma 3.1(i)]{HKM06}.

%

\noindent (ii): We first consider the case where $\phi_2$ is not a vertex of 
$B(\s,\alpha)$. Then there exists some split $S = A|B \in \s -\s'$
such that $\phi_2(A)=\alpha$ and $\phi_2(B) =\frac{\alpha(S)}{2} - \alpha$ where
$0 < \alpha <\frac{\alpha(S)}{2} $. But, by (i), 
$\phi_1(A) \in \{0,\frac{\alpha(S)}{2}\}$
for all $A \in S$. Hence, without loss of generality, $\phi_1(A)=0$ 
and $\phi_2(A) \neq 0$.

Now, suppose that $\phi_2$ is a vertex of $B(\s,\alpha)$. 
Then as $\phi_2 \in B(\s,\alpha) - B_{\s'}(\s,\alpha)$ there is a path
$\phi, \psi_1,\psi_2, \dots , \psi_m=\phi_2$, $m \ge 1$,  in the
1-skeleton of  $B(\s,\alpha)$ where $\phi$ is a cut-vertex of  $B(\s,\alpha)$
contained in $B_{\s'}(\s,\alpha)$ such that $\phi_2$ is in one of
the connected components of $B(\s,\alpha) - \{\phi\}$. Let
$S=A|B$ be the split corresponding to the 1-cell in  $B(\s,\alpha)$ with
end vertices $\psi_{m-1}$ and $\phi_2$. 
By Lemma~\ref{zero-one}(i),
$\phi_1|_S = \phi|_S$ follows.
Using the isomorphism between the 0-cells of $B(\s,\alpha)$ and
the  vertices of the Buneman graph $G(\s,\alpha)$ given in 
\cite[Corollary 3.2]{DHM97a}
let $\phi',\phi_2'$ denote the vertices in $G(\s,\alpha)$
corresponding to the vertices $\phi, \phi_2$, 
respectively.
Without loss of generality, $\phi_2'(S)=B$. Hence $\phi_2(B)=0$ and 
$\phi_2(A)=\alpha(S)/2$ in view of that isomorphism. By the choice of $S$,
$\phi'(S)=A$ follows, and so $\phi(A)=0$ using again that
isomorphism. By the choice of $S$, it follows that $0=\phi(A)=\phi_1(A)$.
\epf

\section{A key result on the blocks of the Buneman complex}\label{sec:key}

In this section we will prove a key result which, for a 
weighted split system $(\s,\alpha)$,  allows us to decide whether 
or not two maps $\phi_1, \phi_2 \in B(\s,\alpha)$,  are contained
within the same block of $B(\s,\alpha)$.

We begin with a lemma concerning cut-vertices of the Buneman complex.

\begin{lemma}\label{same}
Suppose that $(\s,\alpha)$ is a weighted split system on $X$, 
$\phi$ is a cut-vertex of $B(\s,\alpha)$, and $C$ is some
connected component\footnote{Here we are using connected component 
in the topological sense -- see e.\,g.\,\cite[p.103]{S83}} 
of $B(\s,\alpha) - \{\phi\}$.\\
(i)  If $S \in \s(C^{+\phi})$, and 
$\phi_1,\phi_2 \in B(\s,\alpha) - C$,  then $\phi_1|_S=\phi_2|_S$.\\
(ii) If $\phi' \in C$ then there
exists some $S \in \s(C^{+\phi})$ such that $\phi'|_S \neq \phi|_S$.
\end{lemma}
\pf (i): In view of Lemma~\ref{restrict}, it suffices to show that 
$\tilde{\phi_1}|_S=\tilde{\phi_2}|_S$
holds for any pair of vertices $\tilde{\phi_1}$ of $[\phi_1]$ 
and $\tilde{\phi_2}$ of $[\phi_2]$. 
But for any such pair $\tilde{\phi_1}$ and $\tilde{\phi_2}$, we 
have that  $\tilde{\phi_1}$ and $\tilde{\phi_2}$ are connected by a 
shortest path in
the 1-skeleton of $B(\s,\alpha) - C$. As, by  \cite[Corollary 3.2]{DHM97a},
the 1-skeleton of $B(\s,\alpha)$ is
isomorphic to the Buneman graph $G(\s',\alpha)$ of $(\s,\alpha)$ 
and for any two vertices 
$\psi$  and $\psi'$ in $G(\s,\alpha)$ we have that $\psi(S)=\psi(S')$ for
all $S'\in \s$ not induced by an edge on a shortest
path from $\psi$ to $\psi'$, it is straight-forward to check that 
$\tilde{\phi_1}|_S=\tilde{\phi_2}|_S$.

(ii): If $\phi'$ is a vertex in $B(\s,\alpha)$, then we can take any split 
$S \in \s$ induced by some edge on a shortest path in the 
1-skeleton of $B(\s,\alpha)$
between $\phi$ and $\phi'$. So assume that  
$\dim([\phi']) \geq 1$. Let $\phi_1$ denote
a vertex in $[\phi']$ such that the number $m$ of edges on a
path $P$ from $\phi_1$ to $\phi$ in the 1-skeleton of $B(\s,\alpha)$
is as small as possible. If $m=0$ then $\phi_1=\phi$ and, so,
there must exist some 
$S\in \s(\phi')\subseteq \s(C^{+\phi})$
such that $\phi'|_S\not=\phi|_S$ as $\phi'$ is not a vertex 
in $B(\s,\alpha)$ whereas $\phi$ is.

If $m\not=0$ then there must exist some edge $e$ on $P$ such that for the split 
$S_e\in \s$ induced by $e$ we have
$\phi_1|_{S_e}\not=\phi|_{S_e}$. Note that $S_e\in \s(C^{+\phi})$, 
as $\phi_1\in C$. 
Since, by the choice of 
$\phi_1$ we have $\phi'|_{S_e}=\phi_1|_{S_e}$, (ii) follows in
this case as well. 
\epf 

%
%


We now state and prove the main result of this section.
For $\phi_1 \neq \phi_2 \in B(\s,\alpha)$, let
$$
\Delta(\phi_1,\phi_2) = \{S \in \s \,:\, \phi_1|_S  \neq \phi_2|_S\}.
$$ 
In addition, for $\phi \in  B(\s,\alpha)$, note that
$$
\s(\phi) = \{ A|B \in \s \,:\, \phi(A), \phi(B) \neq \emptyset \} 
= \{ S \in \s \,:\, S \subseteq supp(\phi) \}.
$$

\begin{theorem}\label{pairs}
Suppose that $(\s,\alpha)$ is a weighted split system on $X$,
$\phi_1 \neq \phi_2 \in B(\s,\alpha)$, and $\s' \in  C({\mathcal I}(\s))$. 
Then $\{\phi_1,\phi_2\} \subseteq B_{\s' }(\s,\alpha)$ if and only if
$\Delta(\phi_1,\phi_2) \subseteq \s'$.
\end{theorem}
\pf
First note that without loss of generality
we can assume $|C(\mathcal I(\s))| \ge 2$, since otherwise $\s =\s'$ and
$B(\s,\alpha)=B_{\s' }(\s,\alpha)$ 
and so the theorem clearly holds.

Assume for contradiction that
$\Delta\{\phi_1,\phi_2\}\subseteq \s'$ but $\{\phi_1,\phi_2\} \not\subset 
B_{\s' }(\s,\alpha)$. It suffices to consider the  following two cases:

\noindent{\bf Case 1:} There is a cut-vertex $\phi$ of $B(\s,\alpha)$ 
so that $\phi_1$ and $\phi_2$ are both contained in $B(\s,\alpha) -C$
for $C$ the  connected component of $B(\s,\alpha) - \{\phi\}$ with
$B_{\s' }(\s,\alpha) \subseteq C^{+\phi}$. Then since 
$\s' \subseteq \s(C^{+\phi})$, 
Lemma~\ref{same}(i) implies for all
$S \in \s'$ that  $\phi_1|_S=\phi_2|_S$. Hence $\Delta(\phi_1,\phi_2)$ 
is not a subset of $\s'$, a contradiction.

\noindent{\bf Case 2:} There is a cut vertex $\phi$ of $B(\s,\alpha)$ 
so that $\phi_1$ is in a connected component $C$ of $B(\s,\alpha) - \{\phi\}$
with $B_{\s' }(\s,\alpha)$ not a subset of $C^{+\phi}$, 
and $\phi_2$ is contained in $D^{+\phi}$ for $D$ the 
connected component of $B(\s,\alpha) - \{\phi\}$
with $B_{\s' }(\s,\alpha) \subseteq D^{+\phi}$. If $\phi_1\not=\phi$
then,
 by Lemma~\ref{same}(ii), there
exists some $S \in \s(C^{+\phi})$ such that $\phi_1|_S \neq \phi|_S$. 
Moreover, as $S \in \s(C^{+\phi})$ we have $S \not\in \s'$
 as $\s(C^{+\phi}) \cap \s(D^{+\phi})=\emptyset$. Hence,
$S \not\in \Delta(\phi_1,\phi_2)$. But this is impossible
since,  by Lemma~\ref{same}(i), $\phi|_S = \phi_2|_S$, and
so $\phi_1|_S \neq \phi_2|_S$.

So assume $\phi_1=\phi$.
 By Lemma~\ref{same}(ii), it follows that there must exist
some $S\in \s(D^{+\phi})$ such that $\phi_1|_S\not=\phi_2|_S$.
Hence, $S\in \Delta(\phi_1,\phi_2) \subseteq \s'$. But this
is impossible since $\s'\subseteq \s( C^{+\phi_1})$ and 
$\s( C^{+\phi_1})\cap \s( D^{+\phi_1})=\emptyset$. 

Conversely, suppose $\{\phi_1,\phi_2\} \subseteq B_{\s' }(\s,\alpha)$.
Let $S \in \s -\s'$. Then there is some cut-vertex $\phi$ 
of $B(\s,\alpha)$ and some connected component $C$
of $B(\s,\alpha) - \{\phi\}$ with $S \in \s(C^{+\phi})$ and 
$\phi_1,\phi_2 \in B(\s,\alpha)  -C$. Hence by Lemma~\ref{same}(i), 
$\phi_1|_S = \phi_2|_S$, i.e. $S \not\in \Delta(\phi_1,\phi_2)$.
\epf

\section{The $\kappa$ map}\label{sec:kappa}

In this section, we begin to relate properties of
the Buneman complex and the tight-span of a totally split-decomposable metric. 

First, recall that a metric $d$ on $X$ is 
{\em totally split-decomposable} if there 
exists a  weighted weakly compatible split system $(\s,\alpha)$ on $X$ with 
$$
d = d_{\s,\alpha} = \sum_{S \in \s} \alpha(S) \delta_S
$$
where, for any split $S$ of $X$ and all $x,y \in X$, $\delta_S(x,y) =1$ 
if $S(x)\not= S(y)$ and $\delta_S(x,y) =0$ else.
If $d$ is such a metric, then it follows by results in \cite{BD92} that 
if $d=d_{\s',\alpha'}$ for some weakly compatible split system $\s'$ and
some weighting $\alpha'$ on $\s'$, then $\s=\s'$ and $\alpha= \alpha'$.
Thus, in what follows we are interested in determining the
polytopal structure of $T(d_{\s,\alpha})$ for $(\s,\alpha)$ 
some weighted weakly compatible split system.
For this, we will exploit properties of a certain map
$\kappa$ from the Buneman complex $B(\s,\alpha)$
to the tight-span $T(d_{\s,\alpha})$ which is defined as follows.

Given a weighted split system $(\s,\alpha)$ on $X$, 
we define the map
$$
\kappa: \R^{\mathcal U(\s)} \to \R^X \,: \, \phi \mapsto (X \to \R \,:\,x\mapsto d_1(\phi,\phi_x) ),
$$ 
where $d_1$ denotes the metric on $\R^{\mathcal U(\s)}$ defined
by setting, for $\phi, \phi' \in  \R^{\mathcal U(\s)}$, 
$$
d_1(\phi,\phi') = \sum_{A \in {\mathcal U(\s)}} |\phi(A) - \phi'(A)|.
$$
Note that $\kappa(B(\s,\alpha))=T(d_{\s,\alpha})$ if and
only if $\s$ is weakly compatible 
\cite{DHM98}. Moreover, 
in case $\s$ is weakly  compatible, the map $\kappa'$ 
defined by taking any maximal cell $C$ in
$B(\s,\alpha)$ to the cell $[\kappa(\phi)]$, where $\phi$ is any generator
of $C$, is a well-defined map which induces a bijection
between the maximal cells of  $B(\s,\alpha)$ 
and the maximal cells of $T(d_{\s,\alpha})$ \cite[Theorem 6.1]{HKM06}

We now prove a useful observation (which is essentially 
also shown in the proof of  (i) $\Rightarrow$ (ii) in \cite[Theorem 7.1]{DHM02b}).

\begin{theorem}\label{octahedral}
Suppose that $(\s,\alpha)$ is a weighted weakly compatible
split system on $X$, and $\phi_1,\phi_2 \in B(\s,\alpha)$ distinct with 
$\kappa(\phi_1)= \kappa(\phi_2)$. If $\phi = (\phi_1+\phi_2)/2$, 
then $\s(\phi)$ is an octahedral split system.
\end{theorem}
\pf
As $H(\s,\alpha)$ is convex, $\phi \in H(\s,\alpha)$. As $\kappa$ is
linear, $\kappa(\phi) = \kappa((\phi_1+\phi_2)/2) = 
(\kappa(\phi_1)+\kappa(\phi_2))/2$, 
which is in $T(d_{\s,\alpha})$. Therefore $\phi \in B(\s,\alpha)$ 
(as $\kappa$ maps $B(\s,\alpha)$ surjectively onto  $T(d_{\s,\alpha})$
and $B(\s,\alpha) = \kappa^{-1}(T(d_{\s,\alpha}))\cap H(\s,\alpha)$ 
since $\s$ is weakly compatible -- see e.\,g.\,\cite[p. 305]{DHM98}). Put
$\s'=\s - \s(\phi)$, $\s''=\s(\phi)$,
 $\kappa^* = \kappa|_{\mathbb R^{\mathcal U(\s')}}$, and
$\kappa'=\kappa|_{\mathbb R^{\mathcal U(\s'')}}$.

Note that, as $\phi_1,\phi_2$ distinct, 
$\emptyset \neq \Delta(\phi_1,\phi_2) \subseteq \s''$.
So, by Lemma~\ref{incomp}, $\s''$ is incompatible.
Thus $B(\s'',\alpha|_{\s''})= H(\s'',\alpha|_{\s''})$,
by \cite[Proposition 3.3]{DHM97}.

Moreover, $\phi_1'=\phi_1|_{\mathcal U(\s'')},  
\phi_2'=\phi_2|_{\mathcal U(\s'')} \in B(\s'',\alpha|_{\s''})$, 
and $\phi_1' \neq \phi_2'$ since
$\Delta(\phi_1,\phi_2) \subseteq \s''$. But 
these considerations imply
$$
\kappa'(\phi_1')= \kappa(\phi_1) - \kappa^*(\phi_1|_{\mathcal U(\s')})
= \kappa(\phi_2) - \kappa^*(\phi_2|_{\mathcal U(\s')})
=\kappa'(\phi_2')
$$ 
since, by \cite[Lemma 3.1(i)]{HKM06},   
$\phi_1|_{\mathcal U(\s')}=\phi_2|_{\mathcal U(\s')}$ holds.
Thus, $\kappa'$ is not injective.
Hence $\s''$ is not strictly circular by \cite[Proposition 5.1]{DHM02b}.

Since $\s''$ is weakly compatible, $\s''$ is
either strictly circular or octahedral, and so $\s''$ is octahedral.
\epf

Using this last result we now prove two results 
which relate properties of $\kappa$ to the blocks of the
Buneman complex.
The first result shows that 
$\kappa$ is injective when restricted to any 
block of the Buneman complex that does not correspond
to an octahedral split system.

\begin{theorem}\label{consistentblock}
Suppose that $(\s,\alpha)$ is a weighted weakly compatible
split system on $X$ and that $\s' \in C(\mathcal{I}(\s))$ is
not octahedral, then $\kappa|_{B_{\s'}(\s,\alpha)}$ is injective.
\end{theorem}
\pf
Suppose $\phi_1 \neq \phi_2 \in B_{\s'}(\s,\alpha)$ with 
$\kappa(\phi_1)=\kappa(\phi_2)$.
By Theorem~\ref{pairs}, 
$\emptyset\not=\Delta(\phi_1,\phi_2) \subseteq \s'$.
Put $\phi=(\phi_1+\phi_2)/2$. 
By Theorem~\ref{octahedral},  $\s(\phi)$ is
an octahedral split system. 
Since $\s' $ is not octahedral we cannot have
$\s(\phi)\subseteq \s'$, by Corollary~\ref{structure}.
Hence, $\Delta(\phi_1,\phi_2)\subseteq 
\s(\phi)\cap S'=\emptyset$ which is impossible.
\epf

The second result shows that if $\phi_1,\phi_2 \in B(\s,\alpha)$ are distinct
and not in the same block of $B(\s,\alpha)$, then 
their images under $\kappa$ are distinct.

\begin{theorem}\label{distinctmap}
Suppose that $(\s,\alpha)$ is a weighted weakly compatible
split system on $X$.
If $\phi_1,\phi_2 \in B(\s,\alpha)$ are distinct, 
and $\{\phi_1,\phi_2\} \not\subseteq B_{\s'}(\s,\alpha)$,
 for any $\s' \in C(\mathcal{I}(\s))$,
then $\kappa(\phi_1)\neq \kappa(\phi_2)$.
\end{theorem}
\pf
Suppose for contradiction that $\kappa(\phi_1)=\kappa(\phi_2)$. 
Let $\phi=(\phi_1+\phi_2)/2$. By 
Theorem~\ref{octahedral}, $\s(\phi)$ is an octahedral split system.
Since $\s(\phi)\subseteq \s'$ must hold for some $\s'\in  C(\mathcal{I}(\s))$, 
Corollary~\ref{structure} implies $\s'=\s(\phi)$. Thus,
$\Delta(\phi_1,\phi_2)\subseteq \s'$ and so, 
by Theorem~\ref{pairs}, $\{\phi_1,\phi_2\} \subseteq  B_{\s'}(\s,\alpha)$, 
a contradiction.
\epf

\section{Relating blocks in $B(\s,\alpha)$ to blocks in $T(d_{\s,\alpha})$}\label{sec:relate}

In this section we use the map $\kappa$ to provide an
explicit bijection between the blocks of the Buneman complex
and the tight-span for the metric $d_{\s,\alpha}$ associated
to a weighted weakly compatible split system $(\s,\alpha)$.

We begin by proving that $\kappa$ maps the underlying set of any maximal cell
in the Buneman complex onto the underlying set of some maximal cell in the tight-span.

\begin{lemma}\label{maximal}
Suppose that $(\s,\alpha)$ is a weighted weakly compatible
split system on $X$, $[\phi]$ 
is a maximal cell in $B(\s,\alpha)$ with generator 
$\phi \in  B(\s,\alpha)$, 
and $\kappa'([\phi]) = C$, with $C$ a maximal cell in $T(d_{\s,\alpha})$. Then
$\kappa([\phi]) = C$.
\end{lemma}
\pf First note that by the definition of $\kappa'$, $\kappa(\phi)$ is
a generator for $C$, i.e. $C=[\kappa(\phi)]$.

Now, suppose  $\psi \in [\phi]$, then by \cite[Theorem 5.1(iii)]{HKM06}
$\kappa(\psi) \in [\kappa(\phi)] = C$. Hence $\kappa([\phi]) \subseteq C$.

Conversely, suppose $f \in C$. Since $\kappa$ is surjective, there is some $\psi \in B(\s,\alpha)$
with $\kappa(\psi)=f$. But as $\kappa(\psi) \in [\kappa(\phi)]$, by 
\cite[Theorem 5.1(iii)]{HKM06}, $\psi \in [\phi]$. Hence $f \in \kappa([\phi])$, and
so $C \subseteq \kappa([\phi])$.
\epf

We now show that the image under $\kappa'$ of 
any pair of maximal cells in the Buneman complex 
can intersect in at most one point.

\begin{lemma}\label{maxmiss}
Suppose that $(\s,\alpha)$ is a weighted weakly compatible
split system on $X$, $\s' , \s'' \in C({\mathcal I}(\s))$ distinct, and
that $\Omega, \Omega'$ are maximal cells in $B_{\s'}(\s,\alpha)$ and
$B_{\s''}(\s,\alpha)$, respectively. Then 
$|\kappa'(\Omega) \cap \kappa'(\Omega')|\le 1$.
\end{lemma}
\pf
Let $C = \kappa'(\Omega)$ and $C'=\kappa'(\Omega')$ and, for
the purposes of obtaining a contradiction, that $|C \cap C'| > 1$.
Then as $C \cap C'$ is a cell in $T(d_{\s,\alpha})$, $dim(C \cap C') \ge 1$.

Now, suppose that $g,g'$ are distinct generators for $C \cap C'$, which
must exist since  $dim(C \cap C') \ge 1$. Since $\kappa(\Omega)=C$
and  $\kappa(\Omega')=C'$ and, by \cite[Theorem]{DHM98}, 
$\kappa$ is surjective 
there must exist $\phi, \phi' \in \Omega$ distinct with $\kappa(\phi)=g$ and 
$\kappa(\phi')=g'$, and $\phi'', \phi''' \in \Omega'$ distinct 
with $\kappa(\phi'')=g$ and $\kappa(\phi''')=g'$.
 
Note that $|\Omega \cap \Omega'| \le 1$
since $\Omega, \Omega'$ are contained in distinct blocks of $B(\s,\alpha)$.
If $|\Omega \cap \Omega'| = 0$, then we obtain a contradiction
to Theorem~\ref{distinctmap}, since $\kappa(\phi)=g = \kappa(\phi'')$ 
and $\{\phi,\phi''\}$
is not contained in any block of $B(\s,\alpha)$. Moreover, if 
$|\Omega \cap \Omega'| = 1$,
then at least one of $\phi\not=\phi''$ and $\phi'\not=\phi'''$
must hold. Without loss of generality we may assume that
 $\phi\not=\phi''$. Then $\{\phi,\phi''\}$ is not contained in any block
of $B(\s,\alpha)$. By Theorem~\ref{distinctmap},
$g=\kappa(\phi)\not=\kappa(\phi'')=g$ which is also impossible.
\epf

We now show that if a block in the Buneman complex
consists of a single maximal cell, then its image under
$\kappa$ is the underlying set of some block in the tight-span.

\begin{proposition}\label{cellblock}
Suppose that $(\s,\alpha)$ is a weighted weakly compatible
split system on $X$, $\s' \in C({\mathcal I}(\s))$, and the block 
$B_{\s'}(\s,\alpha)$ consists of a single maximal cell $\Omega$ 
in $B(\s,\alpha)$. Then $\kappa'(\Omega)$ is a block in $T(d_{\s,\alpha})$.
\end{proposition}
\pf
First note that $C = \kappa'(\Omega)$ is a maximal cell in 
$T(d_{\s,\alpha})$. Suppose that $C$ is not a block of $T(d_{\s,\alpha})$. Then there is
some cell $C'$ distinct from $C$ in $T(d_{\s,\alpha})$ with $C' \not\subseteq C$
and $C \cap C'$ is a cell with dimension at least 1. Let $C''$ be a 
maximal cell in $T(d_{\s,\alpha})$ containing $C'$. Note that $C'' \neq C$
and $C \cap C' \subseteq C \cap C''$. Let $\Omega'\neq \Omega$ be a maximal
cell in  $B(\s,\alpha)$ with $\kappa'(\Omega')= C''$, which
exists since, by \cite[Theorem 6.1]{HKM06}, $\kappa'$ 
is bijecitve. Note that since $B_{\s'}(\s,\alpha)$
contains a single maximal cell, $\Omega'$ and $ \Omega$ must be contained
in different blocks of $B(\s,\alpha)$. Since  
$|\kappa'(\Omega) \cap \kappa'(\Omega')| > 1$,
this is impossible in view of Lemma~\ref{maxmiss}.
%
\epf

We now extend the previous result, and 
show that the image under $\kappa$ of any 
block of the Buneman complex is the underlying set
of some block in the tight-span.

\begin{theorem}\label{equalblock}
Suppose that $(\s,\alpha)$ is a weighted weakly compatible
split system on $X$ and $\s' \in C({\mathcal I}(\s))$.
Then $\kappa(B_{\s'}(\s,\alpha))$ is 
equal to some block of  $T(d_{\s,\alpha})$.
\end{theorem}
\pf
By Proposition~\ref{cellblock} and Corollary~\ref{structure}, 
it suffices to assume that $\s'$ is consistent and that 
$B_{\s'}(\s,\alpha)$ contains 
at least 2 maximal cells.

We first show that $\kappa(B_{\s'}(\s,\alpha))$ is 
a subset of some block of $T(d_{\s,\alpha})$.
Suppose that $\Omega,\Omega'$ are two distinct maximal cells
in $B_{\s'}(\s,\alpha)$ with $dim(\Omega \cap \Omega') \ge 1$, which 
must clearly exist as $B_{\s'}(\s,\alpha)$ is a block of $B(\s,\alpha)$.
We claim that $dim(\kappa'(\Omega) \cap \kappa'(\Omega')) \ge 1$.

To see that this claim holds, first note that by Lemma~\ref{maximal}
$$
\kappa(\Omega \cap \Omega') 
\subseteq \kappa(\Omega) \cap \kappa(\Omega')
=\kappa'(\Omega) \cap \kappa'(\Omega').
$$
Hence, $\kappa'(\Omega) \cap \kappa'(\Omega') \neq \emptyset$
as $dim(\Omega\cap\Omega')\geq 1$ and so $ \Omega\cap\Omega'\not=\emptyset$.
Suppose $dim(\kappa'(\Omega) \cap \kappa'(\Omega')) < 1$. Then 
$\kappa'(\Omega) \cap \kappa'(\Omega') = \{g\}$, with $g$ a vertex of  
$T(d_{\s,\alpha})$. Thus,
$\kappa(\Omega \cap \Omega') = \{g\}$.
But, since $dim(\Omega \cap \Omega') \ge 1$, there must
exist $\phi,\phi' \in \Omega \cap \Omega' \subseteq B_{\s'}(\s,\alpha)$ 
distinct
with $\kappa(\phi)=\kappa(\phi')=\{g\}$, which contradicts
 Theorem~\ref{consistentblock} and therefore proves the claim.

Containment of $\kappa(B_{\s'}(\s,\alpha))$ in some block of $T(d_{\s,\alpha})$
now follows since, as $B_{\s'}(\s,\alpha)$ is
a block of $B(\s,\alpha)$, we can find an ordering
$\Omega_1,\Omega_2,\dots\Omega_k$, $k\ge 2$, of
the maximal cells of $B_{\s'}(\s,\alpha)$ so that for all $1 < i \le k$, there
exists some $1 \le j < i$ such that $dim(\Omega_i \cap \Omega_j) \ge 1$.
Hence by the claim and Lemma~\ref{maximal} it immediately follows that 
$$
\kappa(B_{\s'}(\s,\alpha)) \subseteq \bigcup_{i=1}^k \kappa(\Omega_i) = 
\bigcup_{i=1}^k \kappa'(\Omega_i)
$$
is contained in some block $\mathcal T$ of $T(d_{\s,\alpha})$.

Now,  suppose for contradiction that 
$\kappa(B_{\s'}(\s,\alpha))$ is a strictly proper subset 
of $\mathcal T$. Then 
there must exist two maximal cells $C,C'$ in $T(d_{\s,\alpha})$
with $dim(C \cap C') \ge 1$, $C \subseteq  \kappa(B_{\s'}(\s,\alpha))$ and
$C' \subseteq  \mathcal T - (\kappa(B_{\s'}(\s,\alpha))\cup (C\cap C'))$. Let 
$\Omega, \Omega'$ be distinct maximal cells in $B(\s,\alpha)$
with $\kappa'(\Omega)=C$ and $\kappa'(\Omega')= C'$. 
Then $\Omega \subseteq B_{\s'}(\s,\alpha)$ 
(since $C \subseteq \kappa(B_{\s'}(\s,\alpha))$),
and $\Omega' \subseteq  B_{\s''}(\s,\alpha)$ with 
$\s'' \in C({\mathcal I}(\s))-\{\s'\}$
(since $C' \not\subseteq  \kappa(B_{\s'}(\s,\alpha))$). But $dim(C \cap C') \ge 1$
and so, in view of the claim, we obtain  a contradiction to Lemma~\ref{maxmiss}.
\epf

Putting our results together we now show that,
by mapping the underlying set of each block of the  Buneman complex
onto the underlying set of some block in the tight-span,  the map 
$\kappa$ in fact induces a bijection from the blocks of the
Buneman complex onto the blocks of the tight-span.
Before proving this we recall some definitions and results from \cite{HKM06}.

First, for any cell $\Omega$ 
in the Buneman complex $B(\s,\alpha)$ 
and any $x\in X$ there exists some map $\gamma^x\in \Omega$
called the {\em gate for $x$ in $\Omega$} such that
$d_1(\phi_x,\psi)=d_1(\phi_x,\gamma^x)+d_1(\gamma^x,\psi)$, 
for all $\psi\in \Omega$. Second, to any $x\in X$ associate the map 
$h_x:X\to \R$ in $T(d_{\s,\alpha})$
given by $ h_x(y)= d_{\s,\alpha}(x,y)$, for all $y\in X$. Then, 
similarly, for any cell $C$ in $T(d_{\s,\alpha}) $ and any $x\in X$
there exists  a (necessarily unique) map $g^x\in C$ called
the {\em gate of $C$ for $x$} such that  
$d_{\infty}(h_x,h)=d_{\infty}(h_x,g^x)+ d_{\infty}(g^x,h)$, for all $h\in C$.

We now prove the aforementioned result. 
 
\begin{theorem}\label{Kmap}
Suppose that $(\s,\alpha)$ is a weighted weakly compatible
split system on $X$.  The map 
$$
K: {\B}(B(\s,\alpha)) \to {\B}(T(d_{\s,\alpha})): 
B_{\s'}(\s,\alpha) \mapsto \kappa(B_{\s'}(\s,\alpha))
$$
is a well-defined bijection. Moreover, 
$K$ induces a bijection 
between the set of blocks $\{B_{\s'}(\s,\alpha) \,:\, \s' \in Oct(\s)\}$
in $B(\s,\alpha)$ and the 
set of cells in $T(d_{\s,\alpha})$ that are rhombic dodecahedra.
In particular, if a cell in $T(d_{\s,\alpha})$ is a rhombic dodecahedron, then 
it must also be a block of $T(d_{\s,\alpha})$.
\end{theorem}
\pf
By Theorem~\ref{equalblock}, the map $K$ is well-defined.

To see that $K$ is surjective, suppose 
$\mathcal T \in  {\B}(T(d_{\s,\alpha}))$.
Let $C$ be some maximal cell in $\mathcal T$. Let $\Omega$ be the maximal 
cell in $B(\s,\alpha)$ with $\kappa'(\Omega) = C$, and 
let $B_{\s'}(\s,\alpha)$, $\s' \in C({\mathcal I}(\s))$,
be the block in  $B(\s,\alpha)$ containing $\Omega$. 
Then Lemma~\ref{maximal}
combined with Theorem~\ref{equalblock} 
implies $K(B_{\s'}(\s,\alpha)) = \mathcal T$.

To see that $K$ is injective, suppose there exist two distinct 
blocks $B_{\s''}(\s,\alpha)$, $B_{\s'}(\s,\alpha)$, $\s', \s'' \in C({\mathcal I}(\s))$,
in  $B(\s,\alpha)$ with $K(B_{\s'}(\s,\alpha))=K(B_{\s''}(\s,\alpha))$.
Then there must exist $f \in K(B_{\s'}(\s,\alpha))$, 
$\phi_1 \in B_{\s'}(\s,\alpha)-B_{\s''}(\s,\alpha)$ and $\phi_2 \in B_{\s''}(\s,\alpha)-B_{\s'}(\s,\alpha)$
with $\kappa(\phi_1) = \kappa(\phi_2)=f$. But this contradicts Theorem~\ref{distinctmap}. Thus, $K$ is a bijection.

Now, let $\mathcal R$ denote the set 
of cells in $T(d_{\s,\alpha})$ that are rhombic dodecahedra.
Suppose $\s' \in Oct(\s)$. 
Then, by Corollary~\ref{structure}, 
$B_{\s'}(\s,\alpha)$ is 
a block of $B(\s,\alpha)$ consisting of a single cell $\Omega$.
Hence, $\Omega$ is a maximal cell of $B(\s,\alpha)$. Let 
$\omega\in B(\s,\alpha)$ denote a generator for $\Omega$. Then
$\s(\omega)\subseteq \s(B_{\s'}(\s,\alpha))=\s'$. Since
$\s(\omega)$ is maximal incompatible and
$\s'$ is incompatible it follows that $\s(\omega)=\s'$.
Hence, $\s(\omega)$ is octahedral. To see that then 
$\kappa'(\Omega)$ must be a rhombic dodecahedron, we
next consider the underlying graphs\footnote{For any finite metric  
space $(Y,d)$  the underlying graph, denoted by $UG(Y,d)$, 
has vertex set $Y$ and edge set consisting 
of those 2-sets $\{x,y\}\subseteq Y$
for which there exists no $z\in Y-\{x,y\}$ such that 
$d(x,y)=d(x,z)+d(z,y)$} 
$UG(\G(\kappa'(\Omega)), d':=d_{\infty}|_{\G(\kappa'(\Omega))})$ -- where
$\G(\kappa'(\Omega))$ is the set of gates of $\kappa'(\Omega)$ -- and 
$UG(\Gamma(\Omega), d_1|_{\Gamma(\Omega)})$ -- 
where $\Gamma(\Omega)$ is the set of gates of $\Omega$.

Note first that by the remark directly following
 the proof of \cite[Claim 2, p.476]{HKM06} $d'$ is cell-decomposable 
as $d_{\s, \alpha}$ is totally split-decomposable. By 
\cite[Theorem 1.1]{HKM04}, $(\G(\kappa'(\Omega)), d')$ is a 
proper antipodal metric space and 
$T(\G(\kappa'(\Omega)), d')$ is polytope isomorphic to $\kappa'(\Omega)$.
Since $\kappa$ induces an isometry between $(\G(\kappa'(\Omega)), d')$
and $(\Gamma(\Omega), d_1|_{\Gamma(\Omega)})$
the proof of \cite[Claim 2, p.476]{HKM06} implies  that 
$UG(\Gamma(\Omega), d_1|_{\Gamma(\Omega)})$
and 
$UG(\G(\kappa'(\Omega)), d')$ are isomorphic graphs.
By \cite[Lemma 3.1]{HKM06},
$d_1(\gamma^x,\gamma^y)=\sum_{S\in \s(\omega)}\alpha(S) \delta_S(x,y)$
for all $x,y\in X$. 
A straight forward check shows that
$UG(\Gamma(\Omega), d_1|_{\Gamma(\Omega)})$ is isomorphic to  $K_{3\times 2}$.
Thus, $UG(\G(\kappa'(\Omega)), d')$ is isomorphic to $K_{3\times 2}$. Hence
$T(\G(\kappa'(\Omega)), d')$ is polytope isomorphic with a rhombic
dodecahedron \cite[Theorem 4.3]{HKM04} and, so $\kappa'(\Omega)$
is polytope isomorphic with a rhombic
dodecahedron.  Therefore 
$K(\{B_{\s'}(\s,\alpha) \,:\, \s' \in Oct(\s)\})  \subseteq \mathcal R$.

Conversely, suppose $R$ is a rhombic 
dodecahedron in $T(d_{\s,\alpha})$. Since $T(d_{\s,\alpha})$
is a polytopal complex and, by \cite[Corollary 7.3]{HKM06}
every cell of $T(d_{\s,\alpha})$ is polytope isomorphic with
either a hypercube
or a rhombic dodecahedron, it follows that $R$ must be a maximal cell
of  $T(d_{\s,\alpha})$. Thus, there exists some maximal cell $\Psi$
of $B(\s,\alpha)$ such that $\kappa'(\Psi)=R$.
Let $\xi\in B(\s,\alpha)$ denote a generator of $\Psi$. 
Assume for contradiction that $\s(\xi)$ is circular.
Then it is straight forward to verify that $UG(\Gamma(\Psi),d|_{\Gamma(\Psi)})$
is a $2m$-cycle where $m=|\s(\xi)|$. Since 
$(\Gamma(\Psi),d_1|_{\Gamma(\Psi)})$ and $(\G(R),d'')$ are isometric metric
spaces, where $d''=d_{\infty}|_{\G(R)}$,
it follows that $UG(\G(R),d'')$ is also a $2m$-cycle. Hence, by
\cite[Theorem 4.2]{HKM04}, $T(\G(R),d'')$ is
an $m$-cube. Thus, $R$ is also an $m$-cube which is impossible. Consequently, 
$\s(\xi)$ must be octahedral. By Corollary~\ref{structure}, 
$\Psi$ is the block $B_{\s(\psi)}(\s,\alpha)$ of $B(\s,\alpha)$.
\epf

\section{The polytopal structure of $T(d_{\s,\alpha})$}\label{sec:blocks}

In this section, we conclude by explaining how to obtain the
polytopal structure of the tight-span of a totally split-decomposable metric
directly from the Buneman complex (and hence from its underlying split system).
We shall do this at the end of the section, but first we need to 
consider the polytopal structure of the blocks in $T(d_{\s,\alpha})$. 

First, note that by Theorem~\ref{Kmap} each maximal cell $R$ in $T(d_{\s,\alpha})$ that 
is a rhombic dodecahedron is also a block of $T(d_{\s,\alpha})$.
Moreover, for such a cell, there exists some $\s' \in Oct(\s)$,
such that $\kappa$ maps the underlying set of the block $B_{\s'}(\s,\alpha)$ 
(which is polytope isomorphic to a 4-cube -- see Figure~\ref{os}) 
onto $R$. In particular,
$\kappa$ restricted to $B_{\s'}(\s,\alpha)$ does not induce a polytope isomorphism
between $B_{\s'}(\s,\alpha)$ and $R$. 
We shall now show that (see Theorem~\ref{theo:structure} below), in contrast to octahedral split systems in $\s$, 
in case $\s' \in C(\mathcal I(\s)) - Oct(\s)$ (i.e $\s'$ is a consistent split system), 
$\kappa$ actually induces a polytopal complex isomorphism between $B_{\s'}(\s,\alpha)$ 
and the block of $T(d_{\s,\alpha})$ that has underlying set $\kappa(B_{\s'}(\s,\alpha))$. 
To prove this we shall require some further terminology.

Let $V$ be a finite-dimensional $\R$-vector space and $P \subseteq V$.
A subset $T \subseteq P$ is an {\em extremal} subset of $P$ if, for
any $u,v \in P$, and any positive real numbers $\gamma, \beta >0$ 
with $\gamma + \beta =1$, the assumption $\gamma u + \beta v \in T$ 
implies $u, v \in T$ (see \cite[p. 51]{DHM00a}). For the following 
proposition we assume $V =  \R^X$.

\begin{proposition}\label{extremal}
Suppose that $(\s,\alpha)$ is a weighted weakly compatible
split system on $X$, and $\s' \in C({\mathcal I}(\s))$. Then
$K(B_{\s'}(\s,\alpha))$ is an extremal subset of $P(d_{\s,\alpha})$.
\end{proposition}
\pf
Let $\mathcal T = K(B_{\s'}(\s,\alpha))= \kappa(B_{\s'}(\s,\alpha))$. 
Suppose $f,g \in P(d_{\s,\alpha})$,
$\gamma, \beta >0$, $\gamma +  \beta =1$, and 
$\gamma f + \beta g \in \mathcal T$.
We need to show that $f, g \in \mathcal T$.

First note that $T(d_{\s,\alpha})$ is an extremal subset of 
$P(d_{\s,\alpha})$. 
Hence $f,g \in T(d_{\s,\alpha})$. Moreover, as $\kappa$
maps $B(\s,\alpha)$ onto $T(d_{\s,\alpha})$ there exist 
$\phi, \phi' \in B(\s,\alpha)$ 
with $\kappa(\phi) =f$ and $\kappa(\phi') =g$. 
Since $\gamma\phi+\beta\phi'\in  H(\s,\alpha)$
clearly holds, it follows in view of 
$B(\s,\alpha) = \kappa^{-1}(T(d_{\s,\alpha})) \cap H(\s,\alpha)$ that
$\gamma \phi + \beta \phi' \in B(\s,\alpha)$. 

Now, as $\kappa$ is linear, we have 
$$
\kappa(\gamma \phi + \beta \phi') =  \gamma \kappa(\phi) + \beta  \kappa(\phi')
=  \gamma f +\beta g \in \mathcal T.
$$
Hence,
by Theorem~\ref{Kmap} we have $\gamma \phi + \beta \phi' 
\in B_{\s'}(\s,\alpha)$.

Assume for contradiction that 
$\{\phi,\phi'\} \not\subseteq  B_{\s'}(\s,\alpha)$. Without loss
of generality $\phi \not\in B_{\s'}(\s,\alpha)$. Put 
$\phi_1 = \gamma \phi + \beta \phi' \in B_{\s'}(\s,\alpha)$ and 
$\phi_2 = \phi \not\in B_{\s'}(\s,\alpha)$.
Then by Lemma~\ref{zero-one}(ii), there exists a split $S\in \s$ 
and some $A\in S$ such that
$$
0 = \phi_1(A) = (\gamma \phi + \beta \phi')(A) = 
\gamma \phi(A) + \beta \phi'(A) 
$$
and $0 \neq \phi_2(A) = \phi(A)$ which is impossible.

Hence $\{\phi,\phi'\} \subseteq  B_{\s'}(\s,\alpha)$. But then
$f = \kappa(\phi), g = \kappa(\phi') \in \mathcal T$, as required.
\epf

To establish the property of the 
map $\kappa$ mentioned at the beginning of this section, 
we shall use the following result from 
\cite{DHM00a} which we recall for the convenience of the reader.

\begin{theorem}\cite[Theorem 1(ii)]{DHM00a}\label{theo:cells}
	Let $V$ and $V'$ be finite dimensional real vector spaces, let 
	$P\subseteq V$ and $P'\subseteq V'$ be convex sets and let $f:V'\to V$ be some affine map with 
	$f(P')\subseteq P$. Moreover let $T$ be some extremal subset of $P$ 
	and let $x\in T':= f^{-1}(T)\cap P'$. If $f$ maps $T'$ bijectively into $T$ 
	then $f$ maps the smallest extremal subset of $P'$ containing $x'$ bijectively onto the 
	smallest extremal subset  of $P$ containing $f(x)$.
\end{theorem}

\begin{theorem}\label{theo:structure}
Suppose that $(\s,\alpha)$ is a weighted weakly compatible
split system on $X$, and $\s'  \in
C(\mathcal I(\s)) - Oct(\s)$. Then
$K(B_{\s'}(\s,\alpha))$ is isomorphic to $B_{\s'}(\s,\alpha)$ 
as polytopal complexes.
\end{theorem}
\pf
Let $V' = \R^{\mathcal U( \s)}$, $V= \R^X$, $P'=H(\s,\alpha)$, $P=P(d_{\s,\alpha})$,
$f = \kappa:V' \to V$, $T = K(B_{\s'}(\s,\alpha))$, and $T'=f^{-1}(T) \cap P'$.
Then it is straight forward to see that $f(P') \subseteq P$ and 
$T'= B_{\s'}(\s,\alpha)$ hold.
Moreover, by Theorem~\ref{extremal}, $T$ is an extremal subset of $P$,
and by Theorems~\ref{Kmap} and \ref{consistentblock}, 
$f$ maps $T'$ bijectively onto $T$. The theorem now follows
by applying Theorem~\ref{theo:cells}.
\epf

We now explain how to determine the polytopal structure of
$T(d_{\s,\alpha})$ from the weighted split system $(\s,\alpha)$.
By the results above, $T(d_{\s,\alpha})$  
has one block for each $\s' \in C(\mathcal I(\s))$, 
and each of these blocks is polytope isomorphic to a rhombic dodecahedron
in case $\s' \in Oct(\s)$ and to the Buneman complex $B(\s',\alpha|_{\s'})$ otherwise
(for example, consider the weighted split-system $(\s,\alpha)$
with $\alpha$ the all-one weight function whose incompatibility 
graph is pictured in Figure~\ref{ic} and whose tight-span $T(d_{\s,\alpha})$
is pictured in Figure~\ref{tspan}).  
Moreover, we can obtain $T(d_{\s,\alpha})$
as a polytopal complex by starting with the Buneman complex $B(\s,\alpha|_{\s})$ 
and replacing each 4-cube in $B(\s,\alpha|_{\s})$ corresponding
to an element in $Oct(\s)$ by a rhombic dodecahedron. This 
is done by first identifying the 6 vertices in the 4-cube
within the Buneman complex
which can be canonically identified as pictured in Figure~\ref{os}. Then
we remove the 4-cube and replace it with  
a rhombic dodecahedron in which the 6 vertices
now correspond to the vertices labelled $2,3,7,4$ and the cut vertices in the 
rhombic dodecahedron forming the block in Figure~\ref{tspan}). 

\subsection*{Acknowledgements}

\noindent The authors thank Andreas Dress for 
some helpful and interesting discussions concerning octahedral split systems.
JK was partially supported by the National Natural Science Foundation of China 
(No. 11471009 and No. 11671376).


\end{document}